\documentclass[10pt,openbib]{article}

\usepackage{amsfonts,amsthm}
\usepackage{amsmath}
\usepackage{amssymb}
\usepackage{graphicx}

\newcommand{\R}{\mathbb{R}}
\newcommand{\C}{\mathbb{C}}
\newcommand{\Z}{\mathbb{Z}}

\newcommand{\dee}{\mathop{\! \, \rm d \!}\nolimits}
\newcommand{\comp}{\, \raisebox{2pt}{$\scriptstyle\circ \, $}}

\newcommand{\setrule}{\, \rule[-4pt]{.5pt}{13pt}\, }

\newcommand{\spann}{\mathop{\rm span}\nolimits}

\newcommand{\ttfrac}[2]{\mbox{$\frac{{\scriptstyle #1}}{{\scriptstyle #2}}$}}

\newcommand{\circdot}{\raisebox{2pt}{\tiny$\, \bigodot$}}

\allowdisplaybreaks

\begin{document}

\title{\textbf{Some metric geometry \\ of the icosahedron}} 
\author{Richard Cushman}
\addtocounter{footnote}{1}
\date{}
\maketitle
\addtocounter{footnote}{1}
\footnotetext{printed: \today}

This paper constructs a Riemann surface associated to the 
icosahedron and discusses the geodesics associated to a flat metric on this surface. 
Because of the icosahedral symmetry, this is 
a distinguished special case of the example treated in \cite{cushman}. 

\section{Geometry of the vertices}

Inscribe the icosahedron $I$ in the $2$-sphere $S^2$ of radius $1$, which we 
identify with complex projective space $\mathbb{CP}$. According Klein \cite[p.60]{klein} its twelve vertices are 
\begin{equation}
\mathrm{Vert} = \{ [0:1], \, [1:0], \, [{\epsilon }^{\nu }(\epsilon + {\epsilon }^4):1], \, 
[{\epsilon }^{\nu }({\epsilon}^2+{\epsilon }^3):1], \, \, \mbox{for $\nu = 0,1,\ldots , 4$} \} , 
\label{eq-zerozero}
\end{equation}
where $\epsilon = {\mathrm{e}}^{2\pi i/5}$. Using stereographic projection from 
the north pole of $S^2$ to its equatorial plane ${\R }^2 = \C $, which is the mapping 
\begin{displaymath}
\mathrm{proj}: \mathbb{CP} \setminus \{ [1:0] \}  \rightarrow \C: [z:1] \mapsto z ,
\end{displaymath}
see \cite{mckean-moll}, the image of $\mathrm{Vert}$ is 
\begin{align*}
\mathrm{vert} & = \{ 0, {\epsilon }^{\nu }(\epsilon + {\epsilon }^4), \, 
{\epsilon }^{\nu }({\epsilon}^2+{\epsilon }^3), \, \, \nu = 0,1,\ldots , 4 \} \\
& = \{ 0, {\epsilon}^{\nu }(2 \, \mathrm{Re}\, {\epsilon}), 
{\epsilon}^{\nu }(2\, \mathrm{Re}\, {\epsilon}^2), \, \, \nu = 0,1,\ldots ,4 \} \\
& = \{ 0, {\epsilon }^{\nu } a , -{\epsilon }^{\nu }b, \, \, \nu =0, 1, \ldots ,4 \} , 
\end{align*}
where $a = 2 \cos \frac{2\pi}{5}$, $b = -2 \cos \frac{4\pi }{5} = 2\cos \frac{\pi }{5}$. Note that $b > a$.   
\par \noindent \hspace{1.4in}\begin{tabular}{l}
\includegraphics[width=180pt]{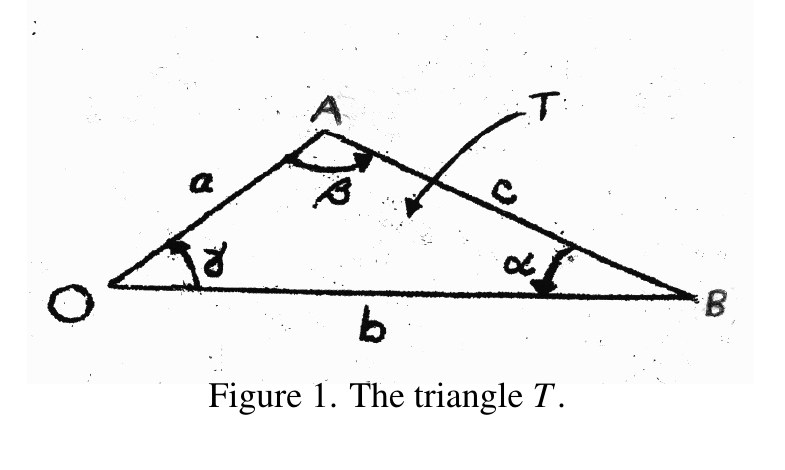}
\end{tabular}

Consider the triangle $T= \triangle BOA$ with the angle $\angle BOA = \gamma = \frac{\pi }{5}$ and sides 
$\overline{OA} = a = 2\cos \frac{2\pi }{5}$,  
$\overline{OB} = b = 2 \cos \frac{\pi }{5}$, and $\overline{AB} =c$, see figure 1. \medskip 

\noindent \textbf{Claim 1.1} $T$ is the rational triangle $(1,2,7)$. \medskip 

\noindent \textbf{Proof.} Applying the law of cosines to $T$ gives 
\begin{align*}
c^2 & = a^2+b^2 -2ab \cos \gamma \\
& = 4 \, {\cos }^2 \ttfrac{2\pi}{5} + 4\, {\cos }^2\ttfrac{\pi }{5} -
8 \, \cos \ttfrac{2\pi }{5}{\cos }^2\ttfrac{\pi }{5} \\
& = 4 \cos \frac{2\pi}{5}\big( \cos \frac{2 \pi}{5} -2{\cos }^2\frac{\pi}{5}\big) +
8 \, {\cos}^2\frac{\pi }{5} - 4 \, {\cos}^2\frac{\pi }{5}\\
& = 4 \, {\sin }^2\ttfrac{\pi }{5}, \, \, \mbox{using $\cos \ttfrac{2\pi}{5} = 2\, {\cos }^2\ttfrac{\pi}{5} -1$ twice.} 
\end{align*}
Thus $c = 2 \sin \frac{\pi }{5}$. Applying the law of sines to the triangle $T$ gives 
\begin{align*}
\sin \beta & = \frac{b}{c} \sin \gamma = 
\frac{2 \cos \frac{\pi }{5} \sin \frac{\pi }{5}}{2\sin \frac{\pi }{5}} = \cos \ttfrac{\pi }{5} >0. 
\end{align*}
So $\beta = \frac{3\pi }{10}$ or $\frac{7\pi}{10}$. If $\beta = \frac{3\pi}{10}$, then 
$T$ is a right triangle, since $\beta + \gamma = \frac{3\pi }{10} + \frac{2\pi }{10} = 
\frac{\pi}{2}$. Hence $\alpha = \beta +\gamma $, which implies $a >b$. This is a contradiction. Hence $\beta = \frac{7\pi }{10}$. So $\alpha = \pi - \gamma - \beta = \pi - \frac{\pi }{5} - \frac{7\pi }{10} = \frac{\pi }{10}$. Thus $T$ is the rational triangle $(1,2,7)$. \hfill $\square $ 
\vspace{-.1in}\par \noindent \hspace{1.1in}\begin{tabular}{l}
\includegraphics[width=150pt]{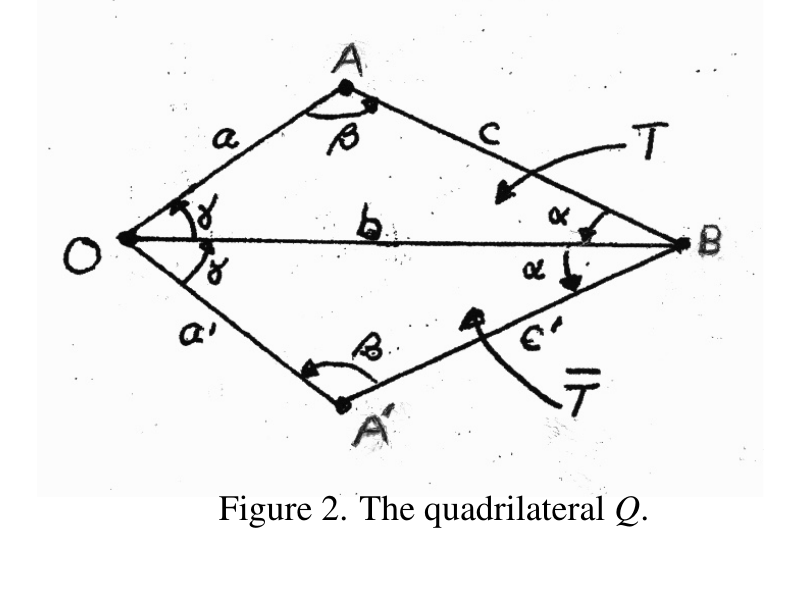} \\
\vspace{-.3in}
\end{tabular}

Let $Q$ be the quadrilateral formed by reflecting the triangle $T$ in its edge 
$\overline{OB}$, see figure 2. Let $K = \bigcup_{0 \le \nu \le 4} R^{\nu}(Q)$, where 
$R: \C \rightarrow \C : z \mapsto \epsilon z = {\mathrm{e}}^{2\pi i/5}z$. 

\par \noindent \hspace{1in}\begin{tabular}{l}
\includegraphics[width=190pt]{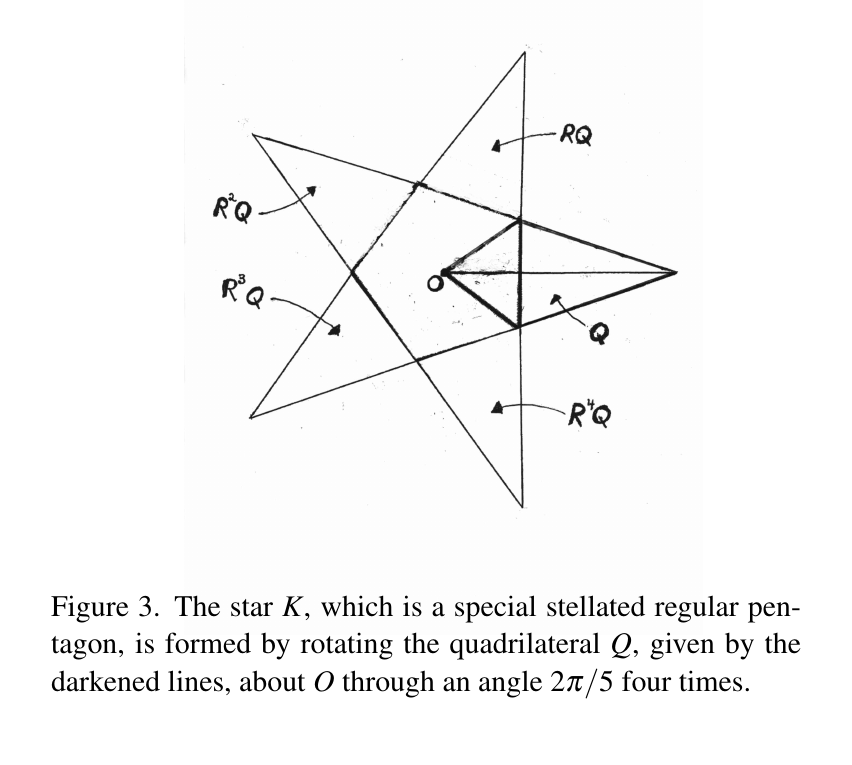}
\end{tabular}

\noindent \textbf{Lemma 1.2} $K$ is a star with center $O$ at $0$, see figure 3. \medskip 

\noindent \textbf{Proof.} Let $\triangle AOB $ be the rational triangle $(2,7,1)$, see 
figure 1.  Form figure 4 \linebreak %

\vspace{-.15in}\noindent \hspace{1.25in}\begin{tabular}{l}
\includegraphics[width=135pt]{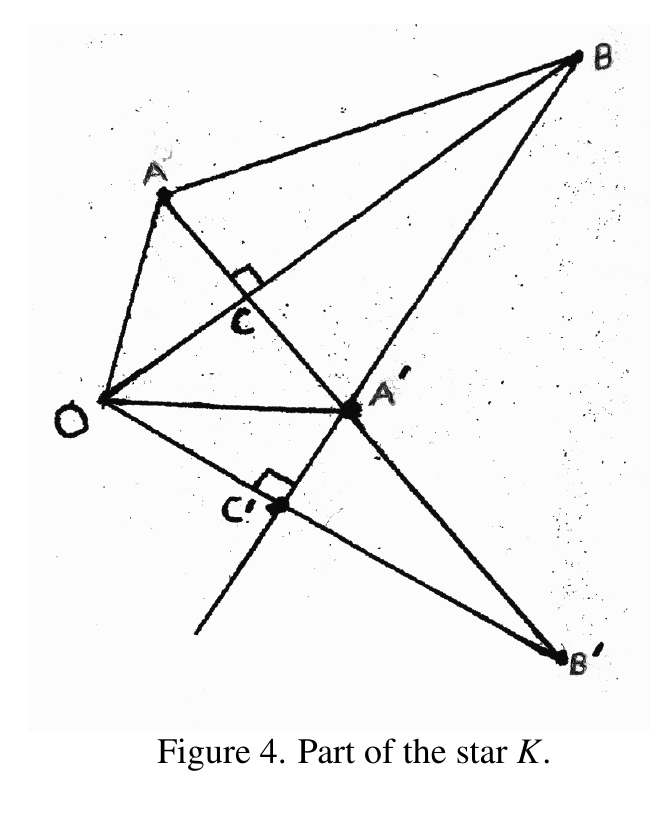}
\vspace{-.1in}
\end{tabular}

\noindent as follows. First reflect the triangle $\triangle OAB$ in the edge 
$\overline{OB}$ to obtain $\triangle OA'B$. Second rotate $\triangle OAB$ about $O$ through an angle 
$2\pi/5$ to obtain $\triangle OA'B'$. \medskip 

We show that the line segments $\overline{AA'}$ and $\overline{A'B'}$ form 
a line segment $\overline{AB'}$. Also the line segments $\overline{BA'}$ and 
$\overline{A'C'}$ form a line segment $\overline{BC'}$. We have the following congruent triangles: 
$\triangle OAB \cong \triangle OA'B \cong \triangle OB'A'$; 
$\triangle BCA' \cong \triangle BCA \cong \triangle B'A'C'$; and 
$\triangle OAC \cong \triangle OA'C \cong \triangle OA'C' $. 
Also we have 
$\angle CAO = \angle CA'O = \angle C'A'O = \frac{3\pi }{10}$; and  
$\angle CAB = \angle CA'B = \angle C'A'B' = \frac{4\pi}{10}$. Thus 
\begin{align*}
\angle CA'B' & = \angle CA'O + \angle OA'C' + \angle C'A'B' 
 = \frac{3\pi}{10} + \frac{3\pi }{10} + \frac{4\pi }{10} = \pi . 
\end{align*}
So $\overline{AB'}$ is a line segment. Also 
\begin{align*}
\angle C'A'B & = \angle C'A'O+ \angle OA'C +  \angle CA'B = 
\frac{3\pi }{10} +\frac{3\pi }{10} + \frac{4\pi }{10} = \pi .
\end{align*}
Hence $\overline{BC'}$ is a line segment. \hfill $\square $ \medskip 

Repeatedly rotating figure 4 about $O$ through $2\pi /5$ gives the star $K$, see figure 3. \medskip

The circumscribed star $\mathcal{K}$ of the star $K$ is formed by constructing the star on the regular 
pentagon formed by the exterior five vertices of $K$. Identifying the exterior 
vertices of the star $\mathcal{K}$, marked by the label $\infty$ in figure 5, and connecting the center of 
$K$, marked by $O$, with all the inner vertices of $K$, and connecting each of the inner vertices of 
$\mathcal{K}$ with an adjacent inner vertex gives figure 5, whose image under the inverse of sterographic projection from the north pole of $S^2$, is the central projection $\widetilde{I}$ of icosahedron $I$ on $S^2$. %

\noindent \hspace{1in}\begin{tabular}{l}
\includegraphics[width=200pt]{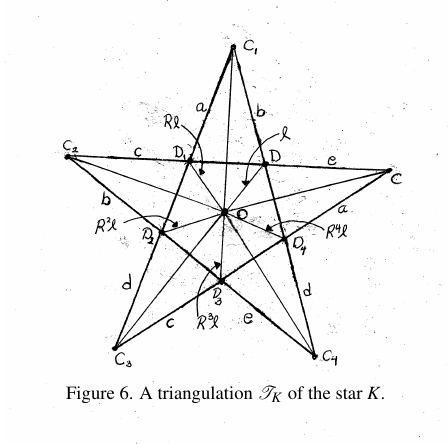}
\end{tabular}

\section{A Riemann surface associated to $K$}

Here we construct an affine Riemann surface $S$ associated 
to the star $K$. \medskip 

The Schwarz-Christoffel mapping associated to the rational triangle $T$ is 
\begin{equation}
F_T: {\C}^{+}\setminus \{ 0, a,b \} \rightarrow T \subseteq \C : 
\xi \mapsto z = F_T(\xi ) = \int^{\xi }_0 k \, \frac{\dee \xi }{\eta }, 
\label{eq-one}
\end{equation}
where $k >0$ and 
\begin{equation}
{\eta }^{10} = {\xi }^{10-2}(\xi -a )^{10-7}(\xi -b)^{10-1} = {\xi }^8(\xi -a)^3(\xi -b)^9 = f(\xi ) , 
\label{eq-two}%
\end{equation}
see \cite{wikipedia}. Here ${\C }^{+}$ is the closed upper half plane 
$\{ \xi \in \C \setrule \, \mathrm{Im}\, \xi \ge 0 \}$. The conformal map $F_T$ is 
a holomorphic diffeomorphism of $\mathrm{int}({\C}^{+}\setminus \{ 0, a,b \})$ 
onto \linebreak
$T^{\ast }= T \setminus \{ O,A,B \}$ with $F_T(0)= O$, $F_T(a) = A$ and $F_T(b) = B$ and is a homeomorphism of $\partial ({\C}^{+}) = \R $ onto  
$\partial T = \overline{OA} \, \cup \, \overline{AB}\, \cup \, \overline{BO}$ with interior 
angles $\angle AOB = \frac{2\pi }{10}$, $\angle OAB = \frac{7\pi}{10}$, and 
$\angle ABO = \frac{\pi }{10}$, see figure 1.\footnote{We determine the constant $k$ in equation (\ref{eq-one}). Let 
$\mathcal{F}(\xi ) = \int^{\xi}_0 \frac{\dee \xi}{\eta }$. Then $\mathcal{F}$ sends 
${\C}^{+}\setminus \{ 0,a,b \}$ onto a triangle $\mathcal{T}$ with interior angles 
$\frac{2\pi }{10}$, $\frac{7\pi}{10}$, and $\frac{\pi }{10}$. The triangles $T$ and 
$\mathcal{T}$ are similar, since they have the same interior angles. Since 
$F_T = k \mathcal{F}$, it follows that $T = k \mathcal{T}$. So 
$k \mathcal{F}(a) = F_T(a) = A = a$, which gives 
\begin{equation}
k = a/\mathcal{F}(a) = a\raisebox{-3pt}{\mbox{\LARGE $/$}} \hspace{-5pt}
\int^{a}_0\frac{\dee \xi }{\eta } . 
\label{eq-three}
\end{equation} }

We now calculate the genus of the compact Riemann surface $\mathcal{S}$ in ${\mathbb{CP}}^2$ \linebreak
associated to the affine Riemann surface $S\subseteq {\C }^2$ defined by equation (\ref{eq-two}), 
see \cite{mckean-moll}. Solving (\ref{eq-two}) for $\eta $ gives
\begin{equation}
\eta = {\epsilon}^{k/2}{\xi }^{1-2/10}(\xi -a)^{1-7/10}(\xi -b)^{1-1/10}, \quad \mbox{$k=0,1, \ldots , 9$} 
\label{eq-threea}
\end{equation}
which defines $\mathcal{S}$ as a branched covering of $\mathbb{CP}$, 
given by 
\begin{equation}
\widetilde{\rho }: \mathcal{S} \subseteq {\mathbb{CP}}^2 \rightarrow \mathbb{CP}  = \C \cup \{ \infty \}: 
([\xi : \eta :1] \mapsto [\xi: 1] 
\label{eq-fivestar}
\end{equation}
with branch values at $0$, $a$, $b$ and possibly $\infty$. At $0$
\begin{subequations}
\begin{align}
\eta & = {\xi }^{1-2/10}(\xi -a)^{1-7/10}(\xi - b)^{1-1/10}  = C_1{\xi }^{1-2/10}(1 + C_2\xi + \cdots ) ; 
\label{eq-threea} 
\end{align}
at $a$ 
\begin{align}
\eta & = (\xi -a)^{1-7/10}(a + (\xi -a))^{1-2/10}(\xi -a -(b-a))^{1-1/10} \notag \\
& = C_1(\xi -a)^{1-7/10}(1 + C_2(\xi -a) + \cdots ); 
\label{eq-threeb}
\end{align}
at $b$ 
\begin{align}
\eta & = (b+ \xi - b)^{1-1/10}(\xi - b + (b-a))^{1-7/10}(\xi -b)^{1-1/10} \notag \\
& = C_1(\xi -b)^{1-1/10}(1 + C_2 (\xi -b ) + \cdots ) ; 
\label{eq-threec}
\end{align} 
at $\infty$
\begin{align}
\eta & = \big( \frac{1}{\xi } \big)^{-(1-2/10)}{\xi }^{1-7/10}(1-\frac{a}{\xi})^{1-7/10}
{\xi }^{1-1/10}(1 - \frac{b}{\xi})^{1-1/10} \notag \\
& = C_1{\xi }^{1-2/10 +1-7/10 +1-1/10}(1+ C_2 {\xi }^{-1} + \cdots ) \notag \\
& = C_1 {\xi }^2(1 + C_2{\xi }^{-1} + \cdots ) .
\label{eq-threed}
\end{align}
\end{subequations}
Therefore $\infty$ is \emph{not} a branch value. The degree of the mapping 
$\widetilde{\rho }$ at the branch value at $0$ is $10/2$; at $a$ is $10/1$ and at 
$b$ is $10/1$. Hence the ramification index $r_0$ at $0$ is 
$2(10/2-1) = 10 -2$; $r_a$ at $a$ is $1(10-1) = 10-1$ and $r_b$ at $b$ is 
$1(10/1-1) = 10 -1$. So the total ramification index $r$ of $\rho $ is 
$r_0 +r_a +r_b = 30 -4$. By the Riemann Hurwitz formula $r = 20 +2g -2$, 
where $g$ is the genus of $\mathcal{S}$. Thus $g =4$. \medskip

We now look more closely at the affine Riemann surface 
$S \subseteq {\C }^2$ (\ref{eq-two}). The singular points of 
$\mathcal{S}$, the closure of $S$ in ${\mathbb{CP}}^2$, is the set 
${\pi }^{-1}(\{ 0, a, b\} ) \cap S$, where $\pi $ is the projection map $\pi : {\C}^2 \rightarrow \C: (\xi ,\eta ) \mapsto \xi $. 
Thus $S_{\mathrm{reg}} = S \setminus \{ (0,0), (a,0), (b,0) \}$ is 
a smooth connected $1$ dimensional complex submanifold of 
${\C }^2 \setminus \{ \eta =0 \}$. The surface $S_{\mathrm{reg}}$ is 
a holomorphic $10$-fold covering of $\C \setminus \{ 0, a,b\} $ given by 
\begin{align}
& {\pi }_{| S_{\mathrm{reg}}}: S_{\mathrm{reg}} \subseteq {\C }^2 \rightarrow \C \setminus \{ 0,a,b\}: 
(\xi , \eta ) \mapsto \xi .
\label{eq-six*}
\end{align}

For $(\xi , \eta ) \in S_{\mathrm{reg}} \subseteq {\C }^2\setminus \{ \eta =0 \}$ the tangent space 
to $S_{\mathrm{reg}}$ at $(\xi ,\eta )$ is
\begin{align*}
T_{(\xi ,\eta )}S_{\mathrm{reg}} &  = \{ \eta \dee \eta - \ttfrac{1}{10} \frac{f^{\prime }(\xi)}{{\eta }^8} \dee \xi =0 \} = 
{\spann }_{\C} \{ \eta \frac{\partial }{\partial \xi } + \ttfrac{1}{10} 
\frac{f^{\prime }(\xi )}{{\eta }^8} \frac{\partial }{\partial \eta } \}  
\end{align*}
with $f'(\xi ) \ne 0$. So the mapping
\clearpage
\begin{align}
&X: S_{\mathrm{reg}} \rightarrow TS_{\mathrm{reg}}:  \notag \\
&\hspace{.25in}(\xi , \eta ) \mapsto 
\big( (\xi ,\eta ), X(\xi ,\eta ) \big) 
 =  \big( (\xi , \eta ), \ttfrac{1}{k}\eta \frac{\partial }{\partial \xi } + 
\ttfrac{1}{10k} \frac{f^{\prime }(\xi )}{{\eta }^8} \frac{\partial }{\partial \eta } \big) 
\label{eq-six**}
\end{align}
is a nowhere vanishing holomorphic vector field on $S_{\mathrm{reg}}$.  \medskip 

Define the mapping 
\begin{equation}
\delta : {S}_{\mathrm{reg}} \subseteq {\C }^2 \rightarrow T^{\ast }\subseteq \C : (\xi , \eta ) \mapsto 
(F_T \comp {\pi }_{|{S}_{\mathrm{reg}}})(\xi , \eta  ) .  
\label{eq-nine*}
\end{equation}
The map $\delta $ is a local holomorphic diffeomorphism of 
${S}_{\mathrm{reg}}$ onto $T^{\ast }$. \medskip 

For $(\xi ,\eta ) \in {S}_{\mathrm{reg}}$ we obtain 
\begin{align}
T_{(\xi ,\eta )}\delta \big( X(\xi ,\eta ) \big) & = 
T_{(\xi ,\eta )} (F_T \comp {\pi }_{|{S}_{\mathrm{reg}}})X(\xi ,\eta ) = 
 T_{\xi }F_T( \frac{1}{k} \eta \frac{\partial }{\partial \xi }) \notag \\
& = \frac{\partial }{\partial z}\rule[-9pt]{.5pt}{18pt}\raisebox{-8pt}{$\scriptscriptstyle \, 
z = \delta (\xi ,\eta ) $}, \hspace{-25pt} 
\label{eq-nine**} 
\end{align} 
since $\dee F_T(\xi ) = k \, \frac{\dee \xi }{\eta } = \dee z$. 
In other words, the conformal map $\delta $ (\ref{eq-nine*}) \emph{straightens} 
the holomorphic vector field $X$ (\ref{eq-six**}) on ${S}_{\mathrm{reg}}$, see \cite{bates-cushman} and \cite{flaschka}. 

\section{Geodesics, symmetry, and billiard motions}

In this section we discuss how the mapping $\delta $ (\ref{eq-nine*}) relates the differential geometry 
of the Riemann surface $S_{\mathrm{reg}}$ to that of $T^{\ast }$. First we 
construct Riemannian metrics $\gamma $ on $\C $ and $\Gamma $ on ${S}_{\mathrm{reg}}$ 
such that the mapping $\delta $ (\ref{eq-nine*}) is an isometry from 
$({S}_{\mathrm{reg}}, \Gamma )$ onto $(T^{\ast }, {\gamma }|_{T^{\ast }})$. We show that a complex time 
rescaling of a real integral curve of the vector field $X$ (\ref{eq-six**}) is 
a geodesic for the metric $\Gamma $. To remove the incompleteness of the geodesic vector field on 
$(T^{\ast }, \gamma |T^{\ast })$ we introduce billiard motions. Next we look at symmetry. 
$K^{\ast}$ is invariant under the 
group $G$, which preserves the metric $\gamma |K^{\ast}$, 
and ${S}_{\mathrm{reg}}$ is invariant under the group $\mathcal{G}$, which 
preserves the metric $\Gamma $. We show that a $G$ invariant extension  
map $\delta $ (\ref{eq-nine*}) intertwines the $\mathcal{G}$ action on $(S_{\mathrm{reg}}, \Gamma )$ with 
the $G$ action on $(K^{\ast }, \gamma |_{K^{\ast }})$. \medskip  

Let $z \in \C$ with $u = \mathrm{Re}\, z$ and $v = \mathrm{Im}\, z$. Then 
\begin{equation}
\gamma = \dee u \circdot \dee u + \dee v \circdot \dee v = 
\dee z \circdot \, \overline{\dee z}
\label{eq-s4fournw}
\end{equation}
is the flat Euclidean metric on $\C $. Consider the flat Riemannian metric $\gamma |_{T^{\ast}}$ on $T^{\ast }$. Pulling back $\gamma |_{T^{\ast }}$ by the mapping $F_T$ (\ref{eq-one}) gives the metric 
$\widetilde{\gamma } = (F_T)^{\ast } \gamma |_{T^{\ast}} =  
k^2{|f(\xi )|}^{-1/5} \dee \xi \circdot \, \overline{\dee \xi } $ on 
$\C \setminus \{ 0, a, b \} $. Pulling $\widetilde{\gamma }$ back by the 
projection mapping $\pi :{\C }^2 \rightarrow \C :(\xi ,\eta ) \mapsto \xi $ gives 
the metric $\widetilde{\Gamma } = {\pi }^{\ast }\widetilde{\gamma } = 
k^2{|f(\xi )|}^{-1/5} \dee \xi \circdot \, \overline{\dee \xi } $ on ${\pi }^{-1}(\C \setminus \{ 0, a, b\} )$. Restricting 
$\widetilde{\Gamma }$ to ${S}_{\mathrm{reg}}$ gives the metric 
$\Gamma = \frac{k}{\eta } \dee \xi \, \circdot \, \frac{k}{\overline{\eta }} \, \overline{\dee \xi } $ 
on ${S}_{\mathrm{reg}}$. Here $\eta = f(\xi )^{1/10}$. \medskip 

\noindent \textbf{Lemma 3.1} The map $\delta = \pi \comp F_T$ is an isometry 
of $({S}_{\mathrm{reg}}, \Gamma )$ onto $(T^{\ast }, {\gamma }|_{T^{\ast }})$. \medskip

\noindent \textbf{Proof.}  We compute. For every $(\xi , \eta ) \in {S}_{\mathrm{reg}}$ we have 
\begin{align}
\Gamma (\xi ,\eta )\big( X(\xi ,\eta ), X(\xi ,\eta ) \big) & = \notag \\
&\hspace{-.75in} = 
\big( \frac{k}{\eta } \dee \xi \big( \ttfrac{1}{k}\eta \frac{\partial }{\partial \xi } + 
\ttfrac{1}{10k} \frac{f^{\prime}(\xi )}{{\eta }^8} \, \frac{\partial }{\partial \eta } \big) \big) \cdot 
 \big( \frac{k}{\overline{\eta }} \overline{\dee \xi }\big(\ttfrac{1}{k} \overline{\eta} \overline{\frac{\partial }{\partial \xi }} + \ttfrac{1}{10k} \frac{\overline{f^{\prime}(\xi )}}{{\overline{\eta }}^8} \, \overline{\frac{\partial }{\partial \eta }} \big) \big)  \notag \\
&\hspace{-.75in} = \frac{k}{\eta } \dee \xi \big( \frac{1}{k}\eta \frac{\partial }{\partial \xi } \big) \cdot  \frac{k}{\overline{\eta }} \overline{\dee \xi } \big(\frac{1}{k} \overline{\eta } \overline{\frac{\partial }{\partial \xi }} \big) =1 = \gamma (z)\big( \frac{\partial }{\partial z}, \frac{\partial }{\partial z} \big)  \notag \\
&\hspace{-.75in} = \gamma (\delta (\xi , \eta ))\big( T_{(\xi ,\eta )}\delta \, X(\xi ,\eta ), 
T_{(\xi ,\eta )}\delta \, X(\xi ,\eta ) \big)  , 
\end{align}
using equation (\ref{eq-nine**}). \hfill $\square $ \medskip

\noindent \textbf{Corollary 3.1A} The metric $\Gamma $ on ${S}_{\mathrm{reg}}$ is flat. \medskip 

\noindent \textbf{Proof.} This follows because the mapping $\delta $ is an 
isometry and the metric ${\gamma }_{T^{\ast }}$ is flat. \hfill $\square $ \medskip 

\noindent \textbf{Corollary 3.1B} The map $\delta $ sends geodesics on 
$({S}_{\mathrm{reg}}, \Gamma )$ to geodesics on 
$(T^{\ast }, {\gamma }|_{T^{\ast }})$ and conversely. \medskip 

\noindent \textbf{Proof.} This follows because the map $\delta $ is an 
isometry. The convese follows because $\delta $ is a local diffeomorphism. 
\hfill $\square $ \medskip 

\noindent \textbf{Corollary 3.1C} The geodesic vector field on 
$({S}_{\mathrm{reg}}, \Gamma )$ is incomplete. \medskip 

\noindent \textbf{Proof.} Since a geodesic motion on $T^{\ast}$ is uniform and rectilinear and $T$ is compact, 
it follows that a geodesic starting at a point in $T^{\ast }$ leaves $T^{\ast }$ in finite time. 
Thus the geodesic vector field on $(T^{\ast }, {\gamma }|_{T^{\ast }})$ is 
incomplete. \hfill $\square $ \medskip 

Since the holomorphic map $\delta $ (\ref{eq-nine*}) straightens the holomorphic 
vector field $X$ on ${S}_{\mathrm{reg}}$, the image under $\delta $ of 
a real time integral curve $\mu $ of $X$ is a real time integral curve of the vector field 
$\frac{\partial }{\partial z}$ on $T^{\ast }$, which is a straight line and hence is a 
geodesic on $(T^{\ast }, {\gamma }|_{T^{\ast }})$. From corollary 3.1B it follows that the 
real time integral curve $\mu $ of $X$ is a geodesic on $({S}_{\mathrm{reg}}, 
\Gamma )$. Let $\alpha \in \C \setminus \{ 0 \}$. For $z_0 \in T^{\ast}$ 
the curve ${\lambda }_{\alpha}: I \subseteq \R \rightarrow T^{\ast } \subseteq \C: t \mapsto \alpha t + z_0$ is a straight line starting at $z_0$ and hence is a geodesic 
on $(T^{\ast }, {\gamma }|_{T^{\ast }})$. The curve ${\lambda }_{\alpha }$ is a real time integral curve of the vector field 
$\alpha \frac{\partial }{\partial z}$ starting at $z_0$. Hence 
${\lambda }_{\alpha }$ is a complex time reparametrization of a real integral curve of 
$\frac{\partial }{\partial z}$. The image of ${\lambda }_{\alpha }$ under the local inverse of the map $\delta $ is a complex time reparametrization of the real time integral curve 
$\mu $ of the vector field $X$ on ${S}_{\mathrm{reg}}$, which is a geodesic on $({S}_{\mathrm{reg}}, \Gamma )$. Thus we have proved \medskip 

\noindent \textbf{Claim 3.2} A complex time integral curve of the holomorphic vector field 
$X$ on ${S}_{\mathrm{reg}}$, which is a reparametrization
of a real time integral curve of $X$, is a geodesic on $({S}_{\mathrm{reg}}, \Gamma )$. \medskip

We now discuss the discrete symmetries of $({S}_{\mathrm{reg}}, \Gamma )$.  
\medskip 

\noindent \textbf{Lemma 3.3} The dihedral group $\mathcal{G} = \langle \mathcal{R}, \mathcal{U} \setrule \, {\mathcal{R}}^5 = e = {\mathcal{U}}^2 \, \, \& \, \, 
\mathcal{U}\mathcal{R} = {\mathcal{R}}^{-1}\mathcal{U} \rangle $, 
where 
\begin{align*}
&\mathcal{R}: {S}_{\mathrm{reg}} \rightarrow {S}_{\mathrm{reg}}: 
(\xi , \eta ) \mapsto ( \xi , \epsilon \, \eta ) \, \, \mathrm{and} \, \, 
\mathcal{U}: {S}_{\mathrm{reg}} \rightarrow {S}_{\mathrm{reg}}: 
(\xi , \eta ) \mapsto (\overline{\xi }, \overline{\eta }) , 
\end{align*}
is a group of isometries of $({S}_{\mathrm{reg}}, \Gamma )$. \medskip 

\noindent \textbf{Proof.} For every $(\xi , \eta ) \in {S}_{\mathrm{reg}}$ we get 
\begin{align}
{\mathcal{R}}^{\ast }\Gamma (\xi ,\eta )\big( X(\xi ,\eta ), X(\xi ,\eta ) \big) & = 
\Gamma \big( \mathcal{R}(\xi , \eta ) \big) \big( T_{(\xi ,\eta )}\mathcal{R}\big( X(\xi ,\eta ) \big) , 
T_{(\xi ,\eta )}\mathcal{R}\big( X(\xi ,\eta ) \big) \big) \notag \\
&\hspace{-1.25in} = \Gamma (\xi, \epsilon \, \eta ) \big( 
 \ttfrac{1}{k}  \eta \frac{\partial }{\partial \xi } + \epsilon \ttfrac{1}{10k} \frac{f^{\prime}(\xi )}{{\eta }^8} \frac{\partial }{\partial \eta },  \ttfrac{1}{k} \eta \frac{\partial }{\partial \xi } + \epsilon \ttfrac{1}{10k} \frac{f^{\prime }(\xi )}{{\eta }^8} 
\frac{\partial }{\partial \eta } \big) \notag \\
& \hspace{-1.25in} = \frac{k^2}{|\eta |^2} \dee \xi \big( \ttfrac{1}{k} \eta \frac{\partial }{\partial \xi } \big)  \cdot 
\overline{\dee \xi } \big( \ttfrac{1}{k} \overline{\eta } \overline{\frac{\partial }{\partial \xi }}  \big)
 = \Gamma (\xi , \eta )\big( X(\xi ,\eta ), X(\xi ,\eta ) \big)  \notag 
\end{align}
and 
\begin{align}
{\mathcal{U}}^{\ast } \Gamma (\xi ,\eta ) \big( X(\xi , \eta ), X(\xi , \eta ) \big) 
& = \Gamma \big(\mathcal{U}(\xi ,\eta ) \big) \big( T_{(\xi ,\eta )}\mathcal{U} 
\big(X(\xi ,\eta ) \big) , T_{(\xi ,\eta )}\mathcal{U} \big( X(\xi ,\eta ) \big) \big) \notag \\
& \hspace{-1.25in} = 
\frac{k^2}{|\eta |^2}  \overline{\dee \xi }(\overline{\ttfrac{1}{k} \eta \frac{\partial }{\partial \xi }} ) \cdot 
 \overline{\overline{\dee \xi }}(\overline{\ttfrac{1}{k} \overline{\eta}} 
\overline{\overline{ \frac{\partial }{\partial \xi }}} ) = 
\Gamma (\xi ,\eta ) \big( X(\xi , \eta ), X(\xi ,\eta ) \big) . 
\tag*{$\square $} 
\end{align}

\noindent \textbf{Lemma 3.4} $G= \langle R, U \setrule \, R^5 = e = U^2 \, \, \& \, \, 
UR = R^{-1}U \rangle $ is a group of isometries of $(K^{\ast} , \gamma |K^{\ast })$. 
Recall that $K^{\ast} = \bigcup^{4}_{v=0}R^{\nu }Q$, where $Q$ is the quadrilateral 
formed by reflecting the triangle $T^{\ast }$ on its edge $\overline{OB}$. \medskip 

\noindent \textbf{Proof.} This follows because $R$ and $U$ are Euclidean motions which 
preserve $K^{\ast }$. \hfill $\square $ \medskip  

\noindent \textbf{Lemma 3.5} $\mathcal{G}$ is the group of 
covering transformations of ${S}_{\mathrm{reg}}$. \medskip

\noindent \textbf{Proof.} The element ${\mathcal{R}}^j{\mathcal{U}}^{\ell}$ of $\mathcal{G}$ sends the 
$k \bmod 10$ sheet ${\eta }_k = {\epsilon}^{k/2} {\xi }^8(\xi -a)^3(\xi -b)^9$ of the covering map 
$\widetilde{\rho}$ (\ref{eq-fivestar}) onto the $\left[ (j -k (\ell \bmod 2))/2 \right] \bmod 10$ sheet. This 
generates all the permutations of the 10 sheets of $\widetilde{\rho}$. \hfill $\square $ \medskip 

Let $\mathcal{D} \subseteq {S}_{\mathrm{reg}}$ be 
the fundamental domain of the action of $\widehat{\mathcal{G}} = \langle \mathcal{R} \setrule \, 
{\mathcal{R}}^5 =e \rangle $ on ${S}_{\mathrm{reg}}$. Define the mapping 
\begin{equation}
{\delta}^{\ast }: S_{\mathrm{reg}} \subseteq {\C}^2 \rightarrow K^{\ast } \subseteq \C : 
(\xi ,\eta ) \mapsto (F_{K^{\ast }} \comp {\pi }_{\mathrm{reg}})(\xi ,\eta ) , 
\label{eq-sixteenstar}
\end{equation}
where for each $\nu =0,1, \ldots , 4$ 
\begin{align*}
&F_{K^{\ast }}: \pi ({\mathcal{R}}^{\nu }(\mathcal{D})) \subseteq \C \setminus \{ 0, a, b \} 
 \rightarrow K^{\ast } \subseteq \C : \xi \mapsto {\epsilon }^{\nu } F_Q(\xi ) 
\end{align*}
and 
\begin{displaymath}
F_Q: \C \setminus \{ 0, a,b \} \rightarrow Q \subseteq \C : \xi \mapsto z = \left\{ 
\begin{array}{rl} F_T(\xi), & \mbox{if $\xi \in {\C}^{+} \setminus \{ 0, a, b  \} $} \\
\rule{0pt}{14pt}\overline{F_T(\overline{\xi})}, & \mbox{if $\overline{\xi} \in \overline{{\C}^{+}\setminus \{ 0,a,b\} }$} 
\end{array}  \right. .
\end{displaymath}
The image of $\mathcal{D}$ under the mapping ${\delta }^{\ast}$ (\ref{eq-sixteenstar}) is the quadrilateral 
$Q$.  \medskip

\noindent \textbf{Lemma 3.6} The mapping ${\delta }^{\ast}$ (\ref{eq-sixteenstar}) intertwines the action $\Phi $ of 
$\mathcal{G}$ on ${S}_{\mathrm{reg}}$ with the action
\begin{equation}
\Psi : G \times K^{\ast } \rightarrow K^{\ast }: (g,z) \mapsto g(z) 
\label{eq-s5one}
\end{equation}
of $G$ on $K^{\ast }$, the star $K$ less its vertices. \medskip 

\noindent \textbf{Proof.} From the definition of the mapping ${\delta }^{\ast}$ we see that for each 
$(\xi ,\eta ) \in {S}_{\mathrm{reg}}$ we have 
${\delta }^{\ast } \big( {\mathcal{R}}^{\nu }(\xi ,\eta ) \big) = 
R^{\nu} \delta (\xi ,\eta )$ for every $0 \le \nu \le 4$. $F_Q(\overline{\xi }) = \overline{F_Q(\xi )}$ by construction and from equation (\ref{eq-six*}) one has $\pi ( \overline{\xi }, \overline{\eta } ) = \overline{\xi }$. From the definition of the mapping ${\delta }^{\ast}$ (\ref{eq-sixteenstar}) we get 
${\delta }^{\ast }(\overline{\xi }, \overline{\eta } ) = \overline{{\delta }^{\ast}(\xi ,\eta )}$ for every 
$(\xi ,\eta ) \in {S}_{\mathrm{reg}}$. In other words, 
${\delta }^{\ast} \big( \mathcal{U}(\xi ,\eta ) \big) = 
U({\delta }^{\ast } (\xi, \eta ))$ for every $(\xi ,\eta ) \in {S}_{\mathrm{reg}}$. 
Hence on ${S}_{\mathrm{reg}}$ we have 
\begin{equation}
{\delta }^{\ast } \comp {\Phi }_g = 
{\Psi }_{\varphi (g)} \comp {\delta }^{\ast } \quad \mbox{for every $g \in \mathcal{G}$.}
\label{eq-s5twostar}
\end{equation} 
The mapping $\varphi : \mathcal{G} \rightarrow G$ sends the generators 
$\mathcal{R}$ and $\mathcal{U}$ of the group $\mathcal{G}$ to the generators $R$ and $U$ of the group $G$, respectively. So it is an isomorphism. \hfill $\square $ 
\medskip

By corollary 3.1C the geodesic vector field on $(K^{\ast}, \gamma |K^{\ast })$ 
is incomplete. We remove this incompleteness by imposing the following 
condition: when a geodesic starting at a point $z_0 \in \mathrm{int}\, K$ meets 
a point $p$ in an edge of $\partial K$, which is not a vertex of $K$, 
then it undergoes a reflection in that edge. Otherwise, it meets $\partial K$ 
as a vertex. Then it reverses its motion. If $z_0 \in \partial K \setminus 
\mathrm{vertex}$, then the geodesic motion reaches a vertex of $K$ in finite 
positive time, after which it reverses its motion. Motions which do not 
reach a vertex of $K$ in finite or infinite time are called \emph{billiard} motions. They are defined for 
all time. Because the metric $\gamma |K^{\ast }$ is invariant under the 
group $\widehat{G} = \langle \mathrm{R}\, \setrule \, R^5 =e \rangle $, all 
billiard motions are invariant under the $\widehat{G}$ action on $K^{\ast }$. Using 
corollary 3.1B and lemma 3.5 we have proved \medskip

\noindent \textbf{Claim 3.7} The image of a $\widehat{G}$ invariant billiard 
motion in $K^{\ast }$ under the local inverse of the developing map 
${\delta }^{\ast }$ (\ref{eq-nine*}) is a $\widehat{\mathcal{G}} = \langle \mathcal{R} \, 
\setrule \, {\mathcal{R}}^5 =e \rangle $ invariant geodesic on 
$({S}_{\mathrm{reg}}, \Gamma )$, which is broken at 
${\{ ({\delta }^{\ast })^{-1}(p_i) \} }_{i \in I}$, where $p_i$ are points on $\partial K$ 
at which the billiard motion undergoes a reflection or a reversal. \bigskip

\vspace{.1in}\noindent {\Large \bf{Appendix} }\bigskip
\ 

\noindent {\Large \bf{A1 A Riemann surface coming from $K$}} \bigskip

We now use the star $K$ to construct a Riemann surface $S^{\dagger}$, 
which is holomorphically diffeomorphic to $S$. \medskip 

Let $R$ be the rotation of $\C $ given by multiplication by 
$\epsilon = {\mathrm{e}}^{2\pi i/5}$ and let $U$ be the reflection given by 
complex conjugation. Set $T =RU$. Then $T$ is a reflection in the closed ray 
$\ell = \{ t{\mathrm{e}}^{2\pi i/5} \, \setrule \, t \in [0,a] \}$ because $T^2$ is the 
identity and $T$ fixes every point on $\ell $. The reflection 
$T_k = R^kTR^{-k} = R^{2k+1}U$, $k=0,1, \ldots , 4$, leaves every point on the closed ray $R^k\ell $ fixed. The 
refections $T_k$, $k=0,1, \ldots , 4$ generate a group. \medskip

Define an equivalence relation on $K \setminus O$ by saying that $x$, $y \in K \setminus O$ are 
\emph{equivalent}, $x \sim y$, if and only if either 1) $x$, $y \in \partial K$ such 
that $y = T_mx$ for some $m\in \{ 0,1, \ldots , 4 \}$ or 2) $x$, $y \in \mathrm{int}\, K \setminus O$ 
and $x=y$. For $x \in K \setminus O$ denote by $[x]$ the set of points in $K\setminus O$ that are equivalent to $x$. 
Let $(K\setminus O)^{\sim }$ be the set of equivalence classes of points in $K\setminus O$ and let 
\begin{equation}
\rho : K \setminus O \rightarrow (K \setminus O)^{\sim}: p \mapsto[p] 
\label{eq-fifteenstar}
\end{equation}
be the map which assigns to $p \in K\setminus O$ its equivalence class $[p]$. Given $K\setminus O$ the topology induced from that of $\C $. Place the quotient topology on $(K\setminus O)^{\sim}$. Then \medskip 

\noindent \textbf{Proposition A1.1}  $(K\setminus O)^{\sim}$ is a connected topological complex manifold of 
dimension $1$ with no boundary and having compact closure in $K^{\sim}$, which is smooth except at points corresponding to the vertices of $K$, where it has a conical singularity. \medskip 

\noindent \textbf{Proof.} To show that $(K^{\ast } \setminus  \mathrm{O}  )^{\sim }$ is a smooth manifold, let $E_{+}$ be an open edge of $K^{\ast }$. For $p_{+} \in E_{+}$ let 
$D_{p_{+}}$ be a disk in $\mathbb{C} $ with center at $p_{+}$, which does not contain a vertex of 
$K =\mathrm{cl}(K^{\ast })$, the closure of $K^{\ast }$ in $\C$. Set $D^{+}_{p_{+}} = K^{\ast } \cap D_{p_{+}}$. 
Let $E_{-}$ be an open edge of $K^{\ast }$, which is equivalent to 
$E_{+}$ via the reflection $T_m$, that is, 
$[\mathrm{cl}(E_{+}), \mathrm{cl}(E_{-}) = T_m(\mathrm{cl}(E_{+}))] \in \mathcal{E}$ is an unordered pair of $T_m$ equivalent edges. Let $p_{-} = T_m(p_{+})$ and set $D^{-}_{p_{-}} 
= T_m(D^{+}_{p_{+}})$. Then $V_{[p]} =  \rho (D^{+}_{p_{+}} \cup D^{-}_{p_{-}}) $ is an open neighborhood of $[p] = [p_{+}] = [p_{-}]$ in $(K^{\ast } \setminus  \mathrm{O} )^{\sim }$, which is a smooth $2$-disk, since the identification mapping $ \rho $ is the identity on $\mathrm{int}\, K^{\ast }$. It follows that 
$(K^{\ast } \setminus  \mathrm{O}  )^{\sim }$ is a smooth $1$-dimensional complex manifold without boundary. 
\par We now handle the vertices of $K$. Let $v_{+}$ be a vertex of $K$ and set $D_{v_{+}} =
\widetilde{D} \cap K$, where $\widetilde{D}$ is a disk in $\mathbb{C} $ 
with center at the vertex $v_{+}= r_0{\mathrm{e}}^{i \pi {\theta }_0}$. The map 
\begin{displaymath}
W_{v_{+}}: D_{+} \subseteq \mathbb{C} \rightarrow D_{v_{+}} \subseteq \mathbb{C}:  
r{\mathrm{e}}^{i \pi \theta } \mapsto |r - r_0| {\mathrm{e}}^{i \pi s (\theta - {\theta }_0)} 
\end{displaymath}
with $r \ge 0$ and $0 \le \theta \le 1$ is a homeomorphism, which sends the wedge with angle 
$\pi $ to the wedge with angle $\pi s$. The latter wedge is formed by the closed edges $E'_{+}$ and $E_{+}$ of 
$K$, which are adjacent at the vertex $v_{+}$ such that ${\mathrm{e}}^{i\pi s}E'_{+} = E_{+}$ with the edge 
$E'_{+}$ being swept out through $\mathrm{int}\, K$ during its rotation to the edge $E_{+}$. Because 
$K$ is a rational regular stellated $5$-gon, the value of $s$ is a rational number for each vertex of 
$K$. Let $E_{-} = T_m(E_{+})$ be an edge of $K$, which is equivalent to $E_{+}$ and set 
$v_{-} = T_m(v_{+})$. Then $v_{-}$ is a vertex of $K$, which is 
the center of the disk $D_{v_{-}} = T_m(D_{v_{+}})$. Set $D_{-} = 
{\overline{D}}_{+}$. Then $D = D_{+} \cup D_{-}$ is a disk in $\mathbb{C} $. The map 
$W: D \rightarrow  \rho (D_{v_{+}} \cup D_{v_{-}} )$, where $W|_{D_{+}} =  
\rho \comp W_{v_{+}}$ and $W|_{D_{-}} =  \rho \comp T_m \comp W_{v_{+}} \comp 
{\mbox{}}^{\overline{\rule{5pt}{0pt}}}$, is a homeomorphism of $D$ into a neighborhood 
$ \rho ( D_{v_{+}} \cup D_{v_{-}})$ of $[v] = [v_{+}] = [v_{-}]$ in 
$(K\setminus O)^{\sim}$. Consequently, the identification space $K^{\sim }$ is a topological manifold.  \hfill $\square $ \medskip 

Let $\widehat{G} = \langle R \, \setrule \, R^5 = e \rangle $. The $\widehat{G}$ 
action 
\begin{displaymath}
\cdot : \widehat{G} \times K \rightarrow K: (g,z) \mapsto g \cdot z = g(z) 
\end{displaymath}
preserves the equivalence relation $\sim $ and is free on $K\setminus O$. Hence 
it induces an action of $\widehat{G}$ on $(K \setminus O)^{\sim}$, which is 
free and proper with orbit map 
\begin{displaymath}
\sigma : (K \setminus O)^{\sim} \rightarrow (K \setminus O)^{\sim}/\widetilde{G} = 
S^{\dagger}: [p] \mapsto [\widehat{G}\cdot p] . 
\end{displaymath}
Thus we have \medskip 

\noindent \textbf{Proposition A1.2} The $\widetilde{G}$ orbit space $S^{\dagger} = 
K^{\sim}/\widehat{G}$ is a compact connected topological complex manifold of dimension $1$. 
The orbit space $S^{\dagger} = (K^{\ast } \setminus \mathrm{O})^{\sim}/\widehat{G}$ is a 
smooth $1$-dimensional complex submanifold of $S^{\dagger}$. 

\noindent \hspace{.8in}\begin{tabular}{l}
\includegraphics[width=250pt]{fig6driver}
\end{tabular}%
 
We now determine the topology of $S^{\dagger}$. The star $K$ less its center $O$ 
has a $\widehat{G}$ invariant triangulation ${\mathcal{T}}_{K\setminus O}$ made up of $10$ open triangles $R^j(\triangle ODC)$ and $R^j(\triangle OD_4C)$ for $j=0, 1, \ldots , 4$; $20$ open edges $R^j(\overline{OC})$, $R^j(\overline{OD})$, $R^j(\overline{CD})$, 
and $R^j(\overline{CD_4})$ for $j = 0,1, \ldots , 4$; and $10$ vertices 
$C_j=R^j(C)$ and $D_j=R^j(D)$ for $j=0,1, \ldots , 4$, see figure 6. Consider the set 
$\mathcal{E} = \{ [E, T_k(E)] , \, k =0, 1, \ldots , 4 \}$ of unordered \emph{pairs} of open 
edges on the boundary of the star $K$. The group $\widehat{G}$ acts on $\mathcal{E}$ by 
$g \cdot [E, T_kE] = [g(E), gT_kg^{-1}\big( g(E) \big)] $, where $g \in \widehat{G}$ 
and $[E, T_kE] \in \mathcal{E}$. Note that $gT_kg^{-1}$ is a reflection in the line $gR^k \ell$ since 
$T_k$ is a reflection in the line $R^k \ell $. Then $\mathcal{E} = \{ a, b,c,d,e \}$, where  
\begin{displaymath}
 \begin{tabular}{ccccc}
$a = [\overline{D_4C}, \overline{D_1C_1} ]$ & \quad & 
$b = [\overline{DC_1}, \overline{D_2C_2} ]$ & \quad & 
$c= [\overline{D_1C_2}, \overline{D_3C_3} ]$ \\
\rule{0pt}{12pt}  $d = [\overline{D_2C_3}, \overline{D_4C_4} ]$ & \quad & 
$e = [\overline{D_3C_4}, \overline{DC} ]$,  & 
\end{tabular}
\end{displaymath}
see figure 6. The $\widehat{G}$ orbit $\mathcal{O}(a)$ of $a \in \mathcal{E}$ is
$\{ R^ja \in \mathcal{E} \setrule \, j =0,1, \ldots , 4 \} = \{ a,c,e \}$ and 
$\mathcal{O}(b) = \{ b,d \}$. Since $\mathcal{E} = \mathcal{O}(a) \cup 
\mathcal{O}(b)$, we have found all the $\widehat{G}$ orbits on $\mathcal{E}$. 
Other $\widehat{G}$ orbits on the edges of ${\mathcal{T}}_{K\setminus O}$ are 
$\mathcal{O}(\overline{OC}) = \{ R^j(\overline{OC}) \setrule \, j=0,1, \ldots , 4 \}$, 
$\mathcal{O}'(\overline{OC}) = \{ R^{(j+1)\! \bmod{5}}(\overline{OC}) \setrule \, 
j =0,1, \ldots , 4 \}$, $\mathcal{O}(\overline{OD})$, 
$\mathcal{O}'(\overline{OD})$, $\mathcal{O}(\overline{CD})$, $\mathcal{O}'(\overline{CD})$, and 
$\mathcal{O}(\overline{CD_4})$, $\mathcal{O}'(\overline{CD_4})$. Looking at the end points of 
the edges in a $\widehat{G}$ orbit in $\mathcal{E}$, we see that 
$\mathcal{O}(D) = \{ D,D_2, D_4, D_1, D_3\}$ is a $\widehat{G}$ orbit of the 
vertex $D$ of the star $K$. Also $\mathcal{O}(C) = \{ C,C_1,C_2,C_3,C_4 \}$ is 
the $\widehat{G}$ orbit of the vertex $C$ of $K$. Since $\mathrm{vert} = 
\mathcal{O}(C) \cup \mathcal{O}(D)$, we have found all the $\widehat{G}$ orbits 
on $\mathrm{vert}$. Hence the triangulation ${\mathcal{T}}_{S^{\dagger}}$ of 
the $\widehat{G}$ orbit space $S^{\dagger}$ is made up of the $\widehat{G}$ 
orbits of $2$ triangles $\triangle ODC$ and $\triangle OD_4C$; the $\widehat{G}$ 
orbits of $2$ edge pairs, namely, $\mathcal{O}(a)$, $\mathcal{O}(b)$; the $\widehat{G}$ orbits of 
the $8$ other edges of ${\mathcal{T}}_{K\setminus O}$, namely,  
$\mathcal{O}(\overline{OC})$, $\mathcal{O}'(\overline{OC})$, 
$\mathcal{O}(\overline{OD})$, $\mathcal{O}'(\overline{OD})$,
$\mathcal{O}(\overline{CD})$, $\mathcal{O}'(\overline{CD})$,
 and $\mathcal{O}(\overline{CD_4})$, $\mathcal{O}'(\overline{CD_4})$; and the $\widehat{G}$ 
orbits of $2$ vertices $\mathcal{O}(C)$ and $\mathcal{O}(D)$. These are the 
faces, edges, and vertices, respectively, of a triangulation 
${\mathcal{T}}_{S^{\dagger}}$ of the $\widehat{G}$ orbit space $S^{\dagger}$. 
Thus the Euler characteristic ${\chi }_{S^{\dagger}}$ of $S^{\dagger}$ is 
$2-10+2 = -6$. Since ${\chi }_{S^{\dagger}} = 2 - 2g$, the genus $g$ of $S^{\dagger}$ is $4$. \bigskip 

\vspace{.1in}\noindent {\Large \bf{A2 An affine model of $S$}}\bigskip

In this section we construct an affine model of the affine Riemann surface 
$S = \{ (\xi ,\eta ) \in {\C}^2 \setrule \, {\eta }^{10} = {\xi }^8(\xi - a)^3(\xi -b)^9 \} $. 
In other words, we find a discrete subgroup $\widehat{\mathfrak{G}}$ of the Euclidean group $\mathrm{E}(2)$ 
of the $2$-plane $\C$, which acts freely and properly on $\C \setminus \mathbb{V}$ and 
has $K^{\ast } \setminus \mathrm{O}$ as its fundamental domain. Here $\mathbb{V}$ is the discrete 
subset of $\C$ formed by translating the center $\mathrm{O}$ and the vertices of the star $K$ by a certain 
finite subgroup $\mathcal{T}$ of translations in $\mathrm{E}(2)$. Identifying certain open edges of 
translates of $K^{\ast } \setminus \mathrm{O}$ by $\mathcal{T}$ gives the identification space 
$(\C \setminus \mathbb{V})^{\sim}$ on which $\widehat{\mathfrak{G}}$ acts freely and properly. We show that 
the $\widehat{\mathfrak{G}}$ orbit space $(\C \setminus \mathbb{V})^{\sim}/\widehat{\mathfrak{G}}$ is holomorphically 
diffeomorphic to $S_{\mathrm{reg}}$.  
\par \noindent \hspace{.6in}\begin{tabular}{l}
\includegraphics[width=250pt]{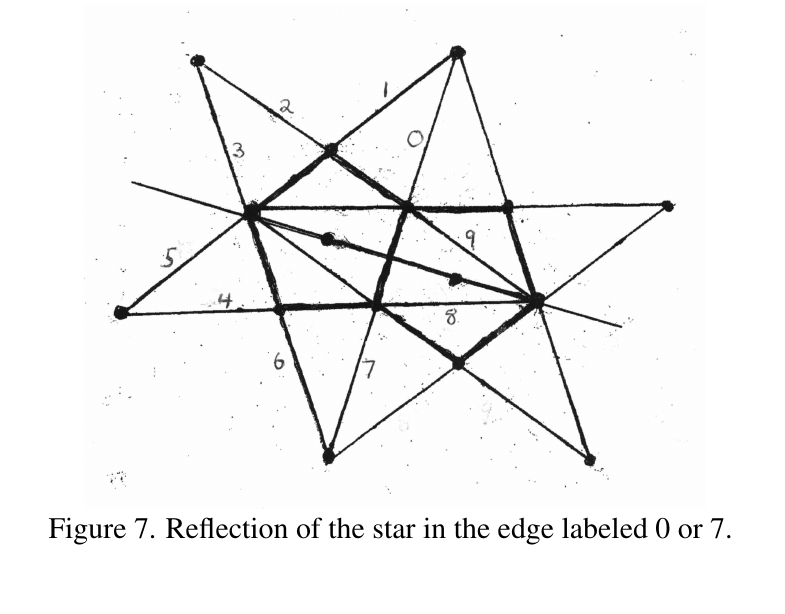}
\vspace{-.1in}
\end{tabular}

First we specify the discrete set $\mathbb{V}$. Give $\C$ its usual orientation. This induces an 
orientation for $\partial K$ and thus a circular ordering of the points of $\partial K$. Label the 
closed edges of $K$ by $\{ 0,1; 2,3; \ldots , 8, 9 \}$. Reflect the star $K$ in one of its edges 
$k_1$. This gives another star $K_{k_1}$ whose center $O_{k_1}$ is the reflection of the center $O$ 
of the star $K$ in its edge $k_1$. Since $K_0=K_7; K_1=K_5; K_2=K_9; K_3=K_6$; and $K_4=K_8$, 
see figure 7, we may describe the above process in a different way. Label the line in the star $K$ by ${\ell }_{k_1}$ 
with $k=0,1, \ldots , 4$, which contains the edge labeled by $k_1$ of the regular pentagon $\mathcal{P}$ 
with center $O$. The star $K_{k_1}$ is obtained by reflecting the star $K$ in the edge $k_1$ of 
$\mathcal{P}$, see figure 7.  The star $K_{k_1}$ is also obtained by translating $K$ so that its 
center is $O_{k_1}$. Repeating this process for each edge of $\mathcal{P}$, and thus each edge of $K$, 
gives $\bigcup^4_{i=0} K_{k_i}$, where $k_i \in \{ 0,1, \ldots , 4\}$, see figure 8.

\par \noindent \hspace{.6in}\begin{tabular}{l}
\includegraphics[width=250pt]{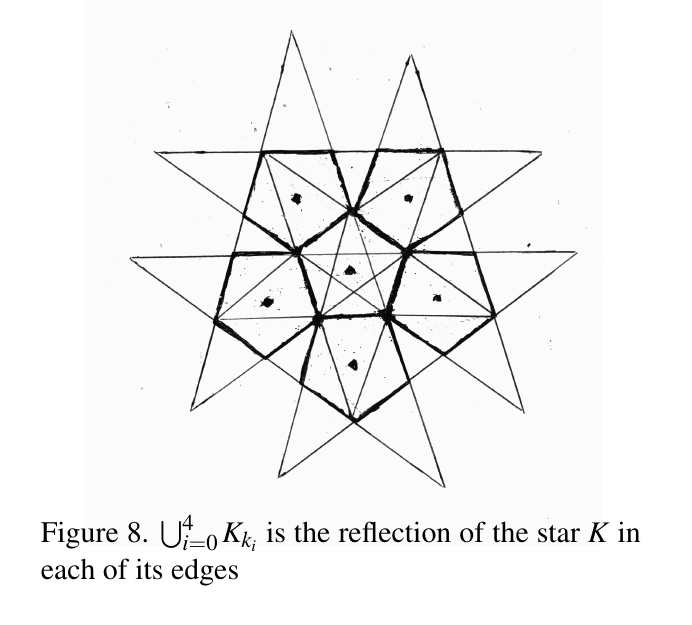}
\end{tabular}
\par \noindent Reflect the star $K_{k_i}$ in the edge $k_j$ of ${\mathcal{P}}_{k_i}$. This gives 
the star $K_{k_ik_j}$, whose center is $O_{k_ik_j}$, which is the reflection of $O_{k_i}$ in the 
edge $k_j$ of ${\mathcal{P}}_{k_i}$. Reflect the regular stellated $5$-gon $K^{\ast} $ in its edge contained in the line 
${\ell }_k$ gives a congruent regular stellated $5$-gon $K^{\ast }_k$ with the center $O$ of $K^{\ast }$ becoming the center $2u_k$ of $K^{\ast }_k$.  The union $\bigcup^4_{i, j =0}{\mathcal{P}}_{k_ik_j}$ is illustrated in figure 9. \medskip

For $k=0,1, \ldots , 4$ let ${\tau }_k$ be the translation 
${\tau }_k: \mathbb{C} \rightarrow \mathbb{C} : z \mapsto z + 2u_k$, 
where $u_k =c \, {\mathrm{e}}^{\frac{(2\pi i}{5})k}$. Here $c$ is the distance 
of the center O of the regular pentagon $\mathcal{P}$ to one of its edges. \medskip

\noindent \textbf{Lemma A2.1} For $k$, $\ell \in \{ 0, 1, \ldots , 4 \} $ we have 
\begin{equation}
{\tau }_{(k+2\ell) \bmod 5} \comp R^{\ell }= R^{\ell } \comp {\tau }_k. 
\label{eq-s4four}
\end{equation}

\noindent \textbf{Proof.} For every $z \in \mathbb{C} $, we have 
\begin{align}
{\tau }_{(k +2\ell) \bmod 5} (z) & = z + 2  u_{(k+2\ell) \bmod 5 } \notag \\
& = z + 2 R^{\ell }u_k  \notag \\
& = R^{\ell} (R^{-{\ell}}z + 2 u_k) = R^{\ell }\comp {\tau }_k (R^{-{\ell}}z). \tag*{$\square $}
\end{align}

\noindent \textbf{Lemma A2.2} The collection of all the centers of the regular stellated $5$-gons, formed by reflecting 
${K^{\ast}}$ in its edges and then reflecting in the edges of the reflected regular stellated $5$-gons et cetera, is 
\begin{align*}
\{ {\tau }^{{\ell }_0}_0 \comp \cdots \comp {\tau }^{{\ell }_{4}}_{4}(0) \in \mathbb{C} \setrule 
({\ell }_0, \ldots , {\ell }_{4}) \in ({\Z }_{\ge 0})^{5} \} =  
\big\{ 2\, \hspace{-20pt} \sum^{\infty}_{\hspace{10pt}{\ell }_0, \ldots , {\ell }_{4} = 0 }
\hspace{-20pt}\big( {\ell }_0 u_0 + \cdots {\ell }_{4}u_{4} \big) \big\}, 
\end{align*}
where for $k=0, 1, \ldots , 4$ we have 
${\tau }^{{\ell }_k}_k = 
\overbrace{{\tau }_k \comp \cdots \comp {\tau }_k}^{{\ell}_k}:\mathbb{C} \rightarrow \mathbb{C}:  z \mapsto z + 
2{\ell }_ju_k$. \medskip

\noindent \textbf{Proof.} For each $k_0 =0, 1, \ldots , 4$ the center of the regular stellated congruent 
$5$-gon $K^{\ast}_{k_0}$ formed by reflecting in an edge of 
$K^{\ast}$ contained in the line ${\ell }_{k_0}$ is ${\tau }_{k_0}(0) = 2u_{k_0}$. Repeating the reflecting process in each edge of $K^{\ast}_{k_0}$ gives congruent regular stellated $5$-gons $K^{\ast}_{k_0k_1}$ with center at 
${\tau}_{k_1}\big( {\tau }_{k_0}(0) \big) = 2(u_{k_1} + u_{k_0})$, where $k_1=0,1, \ldots 4$. Repeating this construction proves the lemma. \hfill $\square $ \medskip 

\par \noindent \hspace{.6in}\begin{tabular}{l}
\includegraphics[width=250pt]{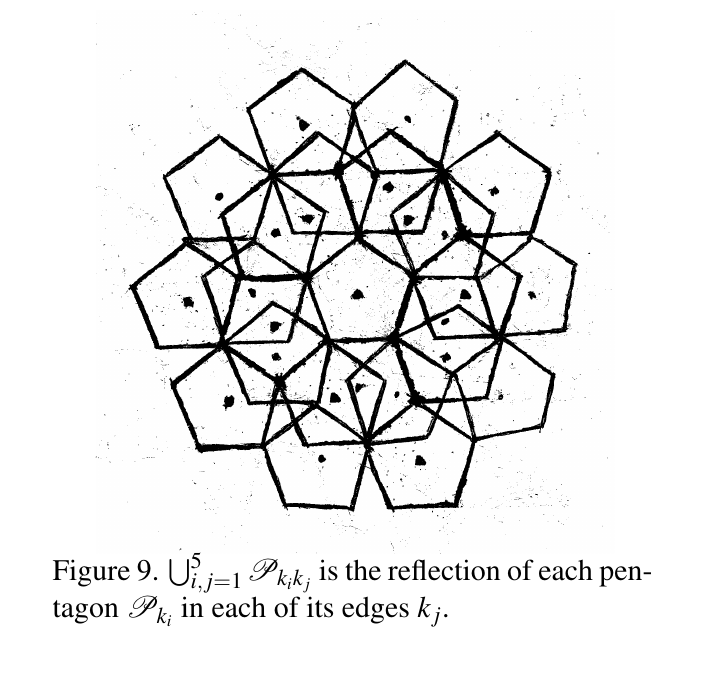}
\end{tabular}

The set $\mathbb{V}$ of vertices of the regular stellated $5$-gon $K^{\ast}$ is 
\begin{displaymath}
\{ V_{2k} = C' {\mathrm{e}}^{2k(\frac{1}{5} \pi \, i )}, \, V_{2k+1} = D'{\mathrm{e}}^{(2k+1)(\frac{1}{5} \pi \, i)} \, \, \mbox{for $0 \le k \le 4$}  \} ,
\end{displaymath}
see figure 3. Clearly the set $\mathbb{V}$ is $G$ invariant. \medskip  

\noindent \textbf{Corollary A2.2A} The set 
\begin{align}
{\mathbb{V}}^{+} & = 
\{ {\tau }^{{\ell }_0}_0 \comp \cdots \comp {\tau }^{{\ell }_{4}}_{4}(V) 
\setrule V \in \mathbb{V} \cup \{ \mathrm{O} \} \, \, \& \, \, 
({\ell }_0, \ldots , {\ell }_{4}) \in ({\Z }_{\ge 0})^{5} \} 
\label{eq-s4threestardot}
\end{align}
is the collection of vertices and centers of the congruent regular stellated 
$5$-gons $K^{\ast }$, $K^{\ast }_{k_1}$, $K^{\ast}_{k_0k_1}, \ldots $. \medskip 

\noindent \textbf{Proof} This follows immediately from lemma A2.2. \hfill $\square $ \medskip 

\noindent \textbf{Corollary A2.2B} The union of $K^{\ast}, K^{\ast}_{k_0}, K^{\ast}_{k_0k_1}, \ldots K^{\ast }_{k_0k_1 \cdots k_{\ell }}, \ldots $, where $\ell \ge 0$, $0 \le j \le \ell $, and $0 \le k_j \le 4$,  covers $\mathbb{C} \setminus {\mathbb{V}}^{+}$, that is, 
\begin{displaymath}
K^{\ast} \cup \bigcup_{\ell \ge 0} \, \bigcup_{0 \le j \le \ell } \, \, \, 
\bigcup_{0 \le k_j\le 4} K^{\ast }_{k_0 k_1 \cdots k_{\ell }} = 
\mathbb{C} \setminus {\mathbb{V}}^{+}. 
\end{displaymath}  

\noindent \textbf{Proof} This follows from the fact that 
$K^{\ast } \cup \bigcup^2_{\ell =0} \bigcup_{0\le j \le \ell} K^{\ast }_{k_0k_1 \cdots k_j}$ covers an 
open neighborhood of $O$, see figure 9, and that $\C \setminus {\mathbb{V}}^{+}$ is obtained by 
repeatedly applying translations ${\tau }_k$ to it. \hfill $\square $ \medskip 

Let $\mathcal{T}$ be the abelian subgroup of the $2$-dimensional Euclidean group 
$\mathrm{E}(2)$ generated by the translations ${\tau }_k$ (\ref{eq-s4three}) for 
$k =0, 1, \ldots 4$. It follows from corollary A2.2B that the regular stellated 
$5$-gon $K^{\ast }$ with its vertices and center removed is the fundamental domain for 
the action of the abelian group $\mathcal{T}$ on $\mathbb{C} \setminus {\mathbb{V}}^{+}$. 
The group $\mathcal{T}$ is isomorphic to the abelian subgroup 
$\mathfrak{T}$ of $(\mathbb{C} , +) $ generated 
by ${\{ 2u_k \} }^{4}_{k=0}$. \medskip

Next we define the group $\mathfrak{G}$ and show that it acts freely, properly, 
and transitively on $\mathbb{C} \setminus {\mathbb{V}}^{+}$. Consider the group $\mathfrak{G} =
G \ltimes \mathfrak{T} \subseteq G \times \mathfrak{T}$, which is 
the semidirect product of the dihedral group $G$, generated by the rotation 
$R$ through $2\pi /5$ and the reflection $U$ subject to the relations 
$R^5 = e = U^2$ and $RU = UR^{-1}$, and the abelian group $\mathfrak{T}$. $\mathfrak{G}$ 
is a discrete subgroup of the Euclidean group $\mathrm{E}(2)$ of the Euclidean motions of the 
$2$-plane. An element $(R^jU^{\ell }, 2u_k)$ of $\mathfrak{G}$ is the affine linear map 
$(R^jU^{\ell }, 2u_k) : \mathbb{C} \rightarrow \mathbb{C} : z \mapsto R^j U^{\ell } z + 2u_k$. 
Multiplication in $\mathfrak{G}$ is defined by 
\begin{displaymath}
(R^jU^{\ell }, 2u_k) \cdot (R^{j'}U^{{\ell}'}, 2u_{k'}) = \big( \mbox{\tiny $\left\{  \begin{array}{rl} 
\hspace{-4pt} R^{j+j'}U^{{\ell}'}, &\hspace{-5pt} \mbox{if $\ell = 0 \bmod 2$} \vspace{3pt}\\
\hspace{-4pt} R^{j-j'}U^{{\ell }'}, & \hspace{-5pt}\mbox{if $\ell = 1 \bmod2$} 
\end{array} \right. $}, (R^jU^{\ell }) (2u_{k'}) + 2u_k \big) ,
\end{displaymath}
which is the composition of the affine linear map $(R^{j'}U^{{\ell }'}, 2u_{k'})$ followed by 
$(R^jU^{\ell }, 2u_k)$. The mappings $G \rightarrow \mathfrak{G}: R^j \mapsto (R^jU^{\ell },0)$ and $\mathfrak{T} \rightarrow \mathfrak{G}: 2u_k \mapsto (e, 2u_k)$ are injective, which allows us to identify the groups $G$ and $\mathfrak{T}$ with their image in $\mathfrak{G}$. Using 
the definition of multiplication we may write an element $(R^jU^{\ell }, 2u_k)$ of $\mathfrak{G}$ as 
$(e, 2u_k) \cdot (R^jU^{\ell },0)$. So 
\begin{displaymath}
(e, 2u_{(j +2k) \bmod 5}) \cdot (R^kU^{\ell },0) = (R^kU^{\ell }, 2u_{(j +2k)\bmod 5}) , 
\end{displaymath} 
The group $\mathfrak{G}$ acts on $\mathbb{C} $ as 
$\mathrm{E}(2)$ does, namely, by affine linear orthogonal mappings. 
Denote this action by  
\begin{displaymath}
\psi : \mathfrak{G} \times \mathbb{C} \rightarrow \mathbb{C} : ((g, \tau ), z ) \mapsto 
\tau (g(z)). 
\end{displaymath} 

\noindent \textbf{Lemma A2.3} The set ${\mathbb{V}}^{+}$ (\ref{eq-s4threestardot}) is 
invariant under the $\mathfrak{G}$ action.  \medskip 

\noindent \textbf{Proof} Let $v \in {\mathbb{V}}^{+}$. Then for some 
$({\ell }'_0, \ldots , {\ell }'_4) \in {\Z }^5_{\ge 0}$ and some 
$w \in \mathbb{V} \cup \{ \mathrm{O} \} $
\begin{displaymath}
v = {\tau }^{{\ell }'_0}_0 \comp \cdots \comp {\tau }^{{\ell }'_4}_{4}(w) = 
{\psi }_{(e, 2u' )}(w) ,
\end{displaymath}
where $u' = \sum^{4}_{k =0} {\ell }'_k u_k$. For $(R^jU^{\ell }, 2u) \in \mathfrak{G}$ 
with $j =0,1, \ldots , 4$ and $\ell =0,1$ we 
have  
\begin{align}
{\psi }_{(R^jU^{\ell }, 2u)}v & = {\psi }_{(R^jU^{\ell }, 2u) } \comp {\psi }_{(e, 2u')} (w) 
= {\psi }_{(R^jU^{\ell },2u) \cdot (e, 2u')}(w) \notag \\
& = {\psi }_{(R^jU^{\ell }, R^j U^{\ell }(2u') +2u)}(w) = 
{\psi }_{(e, 2(R^jU^{\ell }u' +u)) \cdot (R^jU^{\ell },0)}(w) \notag \\ 
& = {\psi }_{(e, 2(R^jU^{\ell }u'+u))} \big( {\psi }_{(R^jU^{\ell },0)} (w) \big)  = 
{\psi }_{(e, 2(R^jU^{\ell }u'+u))}(w'),  
\label{eq-s4seven}
\end{align}
where $w' = {\psi }_{(R^jU^{\ell }, 0)}(w) = R^jU^{\ell }(w) \in 
\mathbb{V} \cup \{ \mathrm{O} \} $. If $\ell =0$, then  
\begin{align}
R^ju' & = R^j (\sum^{4}_{k=0} {\ell }'_k u_k ) = 
\sum^{4}_{k=0} {\ell }'_kR^j (u_k) = 
\sum^{4}_{k=0} {\ell }'_k  u_{(k+ 2j) \bmod 5 };  \notag
\end{align}
while if $\ell =1$, then 
\begin{align}
R^jU(u') & = \sum^{4}_{k=0} {\ell }'_k R^j(U(u_k)) = 
\sum^{4}_{k=0} {\ell }'_kR^j (u_{k'(k)})  
= \sum^{4}_{k=0} {\ell }'_k  u_{(k'(k)+2j)\bmod 5 } . \notag
\end{align}
Here $k'(k) =${\tiny $\left\{ \begin{array}{cl} \hspace{-5pt}2n-k +1, & 
\hspace{-8pt} \mbox{if $k$ is even} \\
\hspace{-5pt} 2n-k-1, & \hspace{-8pt} \mbox{if $k$ is odd}, \end{array} \right. $} \hspace{-5pt}. 
So $(e, 2(R^jU^{\ell }u' +u)) \in \mathfrak{T}$, which implies 
${\psi }_{(e, 2(R^jU^{\ell }u'+u))}(w') \in 
{\mathbb{V}}^{+}$, as desired. \hfill $\square $ \medskip 

\noindent \textbf{Lemma A2.4} The action of $\mathfrak{G}$ on 
$\mathbb{C} \setminus {\mathbb{V}}^{+}$ is free.   \medskip

\noindent \textbf{Proof} Suppose that for some $v \in \mathbb{C} \setminus {\mathbb{V}}^{+}$ and 
some $(R^jU^{\ell }, 2u) \in \mathfrak{G}$ we have $v = {\psi }_{(R^jU^{\ell },2u)}(v)$. Then $v$ lies in some 
$K^{\ast }_{k_0k_1\cdots k_{\ell }} $. So for some 
$v' \in K^{\ast}$ we have 
\begin{align}
v & = {\tau }^{{\ell }'_0}_0 \comp \cdots {\tau }^{{\ell }'_{4}}_{4}(v') 
= {\psi }_{(e, 2u')}(v'),  \notag 
\end{align} 
for some $u' = \sum^{4}_{j=0} {\ell }'_j u_j$ where $({\ell }'_0, \ldots , {\ell }'_{4}) \in ({\Z }_{\ge 0})^5$. Thus, 
\begin{displaymath}
{\psi }_{(e, 2u')}(v') = {\psi }_{(R^jU^{\ell }, 2u) \cdot (e,2u')}(v') = 
{\psi }_{(R^jU^{\ell }, 2R^jU^{\ell }u' +2u)}(v'). 
\end{displaymath}
This implies $R^j U^{\ell }=e$, that is, $j=\ell =0$. So $2u' =2R^jU^{\ell}u' +2u = 2u' +2u$, that is, $u =0$. Hence 
$(R^jU^{\ell },u) = (e,0)$, which is the identity element of $\mathfrak{G}$. \hfill $\square $ \medskip 

\noindent \textbf{Lemma A2.5} The action of $\mathcal{T}$ (and hence 
$\mathfrak{G}$) on $\mathbb{C} \setminus {\mathbb{V}}^{+}$ is transitive. \medskip 

\noindent \textbf{Proof} Let $K^{\ast }_{k_0 \cdots k_{\ell }}$ and 
$K^{\ast }_{k'_0 \cdots k'_{{\ell }'}}$ lie in 
\begin{displaymath}
\mathbb{C} \setminus {\mathbb{V}}^{+} = 
K^{\ast } \cup \bigcup_{\ell \ge 0} \, \bigcup_{0 \le j \le \ell } \, \, 
\bigcup_{0 \le k_j\le 4} K^{\ast }_{k_0 k_1 \cdots k_{\ell }}. 
\end{displaymath}
Since $K^{\ast }_{k_0 \cdots k_{\ell }} = 
{\tau }_{k_{\ell }} \comp \cdots \comp {\tau }_{k_0}(K^{\ast })$ and 
$K^{\ast }_{k'_0 \cdots k'_{{\ell }'}} = 
{\tau }_{k'_{{\ell }'}} \comp \cdots \comp {\tau }_{k'_0}(K^{\ast })$, it follows that 
$({\tau }_{k'_{{\ell }'}} \comp \cdots \comp {\tau }_{k'_0}) \comp 
({\tau }_{k_{\ell }} \comp \cdots \comp {\tau }_{k_0})^{-1}(K^{\ast }_{k_0 \cdots k_{\ell }}) =
K^{\ast }_{k'_0 \cdots k'_{{\ell }'}}$. \hfill $\square $ \medskip   

\noindent The action of $\mathfrak{G}$ on $\mathbb{C} \setminus {\mathbb{V}}^{+}$ is proper because 
$\mathfrak{G}$ is a discrete subgroup of $\mathrm{E}(2)$ with no accumulation points.  \medskip  

We now define an edge of $\mathbb{C} \setminus {\mathbb{V}}^{+}$ and what it means for an 
unordered pair of edges to be equivalent. We show that the group $\mathfrak{G}$ acts freely and 
properly on the identification space of equivalent edges.  

Let $E$ be an open edge of $K^{\ast }$. Since  
$E_{k_0 \cdots k_{\ell }} = {\tau }_{k_0} \cdots {\tau }_{k_{\ell }}(E) \in 
K^{\ast }_{k_0 \cdots k_{\ell }}$, it follows that $E_{k_0 \cdots k_{\ell }}$ is 
an open edge of $K^{\ast }_{k_0 \cdots k_{\ell }}$. Let 
\begin{displaymath}
\mathfrak{E} = \{ E_{k_0\cdots k_{\ell }} \setrule 
\ell \ge 0, \, \, 0 \le j \le \ell  \, \, \& \, \, 0 \le k_j \le 4 \}  
\end{displaymath}
be the set of open edges of $\mathbb{C} \setminus {\mathbb{V}}^{+}$. Since ${\tau }_{k_{\ell }} \comp \cdots 
\comp {\tau }_{k_0}(0)$ is the center of $K^{\ast }_{k_0 \cdots k_{\ell}}$, the element  
$g^{\ast }= (e, {\tau }_{k_{\ell }} \comp \cdots \comp {\tau }_{k_0}) \cdot 
(g, ({\tau }_{k_{\ell }} \comp \cdots \comp {\tau }_{k_0})^{-1})$ of $\mathfrak{G}$ is a 
rotation-reflection of $K^{\ast }_{k_0 \cdots k_{\ell}}$, which sends an edge of 
$K^{\ast }_{k_0 \cdots k_{\ell }}$ to another edge of $g^{\ast }\cdot K^{\ast }_{k_0 \cdots k_{\ell}}$. 
Here $\cdot $ is the usual action of $\mathfrak{G}$ on $\C$. Thus, $\mathfrak{G}$ sends $\mathfrak{E}$ into itself. 
Let ${\mathfrak{E}}_{k_0 \cdots k_{\ell }}$ be the set of unordered pairs 
$[E_{k_0\cdots k_{\ell }}, E'_{k_0\cdots k_{\ell }}]$ of 
equivalent open edges of $K^{\ast }_{k_0 \cdots k_{\ell }}$. In other words, we have
$E_{k_0\cdots k_{\ell }} \cap E'_{k_0\cdots k_{\ell }} = \varnothing $, so the open edges 
$E_{k_0 \cdots k_{\ell }} = {\tau }_{k_0} \cdots {\tau }_{k_{\ell }}(E)$ and 
$E'_{k_0 \cdots k_{\ell }} = {\tau }_{k_0} \cdots {\tau }_{k_{\ell }}(E')$ of 
$K^{\ast }_{k_0 \cdots k_{\ell }}$ are not adjacent, which implies that the 
open edges $E$ and $E'$ of $K^{\ast }$ are not adjacent. Moreover, for some generator 
$T_m$ of the group of reflections with $m \in \{0,1, \ldots , 4 \}$ we have $E' = T_mE$. So
\begin{displaymath}
E'_{k_0 \cdots k_{\ell }} = 
({\tau }_{k_0} \comp \cdots \comp {\tau }_{k_0}) \cdot \big( T_m, ({\tau }_{k_{\ell}} \comp \cdots \comp 
{\tau }_{k_0})^{-1}(E_{k_0 \cdots k_{\ell }} ) \big) . 
\end{displaymath}
Let $\mathcal{E}$ be the 
set of unordered pairs of equivalent edges of $\mathbb{C} \setminus {\mathbb{V}}^{+}$. 
Define an action \raisebox{2pt}{\tiny $\bigstar $} of $\mathfrak{G}$ on $\mathcal{E}$ by 
\begin{align}
(g, \tau ) \raisebox{1pt}{\tiny $\bigstar$} [ E_{k_0 \cdots k_{\ell }}, E'_{k_0 \cdots k_{\ell }} ] & =  
\big( [(\tau  \comp {\tau }' )( g({\tau }')^{-1}(E_{k_0 \cdots k_{\ell }})), 
(\tau  \comp {\tau }' )(g (({\tau }')^{-1}(E'_{k_0 \cdots k_{\ell }})) ] \big) ,  \notag 
\end{align}
where ${\tau }' = {\tau }_{k_{\ell }} \comp \cdots \comp {\tau }_{k_0}$.   

Define a relation $\sim $ on $\mathbb{C} \setminus {\mathbb{V}}^{+}$ as follows. We say that 
$x$ and $y \in \mathbb{C} \setminus {\mathbb{V}}^{+}$ are related, $x \sim y$, if 
1) $x \in F = \tau (E) \in {\mathfrak{E}}$ and $y \in F' = \tau (E') \in {\mathfrak{E}}$ for 
some $\tau \in \mathcal{T}$ and some edges $E$ and $E'$ of $K$ 
such that $[F, F'] = [\tau (E), \tau (E')] =\tau \raisebox{1pt}{\tiny $\bigstar$} [E,E'] \in \mathcal{E}$, where $[E, E'] \in 
{\mathcal{E}}$ with $E' = T_m (E)$  
and $y = \tau \big( T_m({\tau }^{-1}(x)) \big)$ for some $m=0,1 \ldots , 4$, or 2) $x$, $y \in 
\big( \mathbb{C} \setminus {\mathbb{V}}^{+} \big) \setminus \mathfrak{E}$ and $x =y$. 
Then $\sim $ is an equivalence relation on $\mathbb{C} \setminus {\mathbb{V}}^{+}$. Let 
$(\mathbb{C} \setminus {\mathbb{V}}^{+})^{\sim }$ be the set of equivalence classes and 
let $\Pi $ be the map 
$\Pi : \mathbb{C} \setminus {\mathbb{V}}^{+} \rightarrow (\mathbb{C} \setminus {\mathbb{V}}^{+})^{\sim }: 
p \mapsto [p]$, which assigns to every $p \in \mathbb{C} \setminus {\mathbb{V}}^{+}$ the 
equivalence class $[p]$ containing $p$.  \medskip 

\noindent \textbf{Lemma A2.6} $\Pi |_{K^{\ast }}$ is the map $\rho $ (\ref{eq-fifteenstar}).  \medskip

\noindent \textbf{Proof} This follows immediately from the definition of the maps 
$\Pi $ and $\rho $. \hfill $\square $ \medskip 

\noindent \textbf{Lemma A2.7} The action \raisebox{1pt}{\tiny $\bigstar$} of $\mathfrak{G}$ on $\mathbb{C}$, 
restricted to $\mathbb{C} \setminus {\mathbb{V}}^{+}$, is compatible with the equivalence relation 
$\sim $, that is, if $x$, $y \in \mathbb{C} \setminus \mathbb{V}$ and $x \sim y$, then 
$(g, \tau )\cdot (x) \sim (g, \tau )\cdot (y)$ for every $(g, \tau ) \in \mathfrak{G}$. \medskip 

\noindent \textbf{Proof} Suppose that $x \in F = {\tau }'(E)$, where ${\tau }' \in 
\mathcal{T}$. Then $y \in F' = {\tau }'(E')$, since $x \sim y$. So for some 
$T_m$ with $m \in \{ 0,1, \ldots ,4\}$, we have $({\tau }')^{-1}(y) = 
T_m({{\tau }'}^{-1}(x))$ and $F' = T_mE'$. For $(g, \tau ) \in \mathfrak{G}$ then 
$(g, \tau )\cdot (y) \in (g, \tau )\cdot F'$. However,  $(g, \tau )\cdot (x) \in (g, \tau ) \cdot F$ and 
$[ (g, \tau ) \cdot F , (g, \tau ) \cdot F' ]  = (g , \tau )\raisebox{1pt}{\tiny $\bigstar$} [F,F']$. 
Hence $(g, \tau )\cdot (x) \sim (g, \tau )\cdot (y)$.  \hfill $\square $ \medskip 

\noindent \textbf{Lemma A2.8} The action \raisebox{1pt}{\tiny $\bigstar$} of $\mathfrak{G}$ on 
$(\mathbb{C} \setminus {\mathbb{V}}^{+})^{\sim }$ is free and proper.  \medskip 

\noindent \textbf{Proof} 
The following argument shows that the action \raisebox{1pt}{\tiny $\bigstar$} is free. 
An element of $\mathfrak{G}$, which lies in the isotropy group 
${\mathfrak{G}}_{[F,F']}$ for $[F,F'] \in \mathcal{E}$, interchanges the edge $F$ with the equivalent 
edge $F'$ and thus fixes the equivalence class $[p]$ for every $p \in F$. Hence 
the $\mathfrak{G}$ action \raisebox{1pt}{\tiny $\bigstar$} on 
$(\mathbb{C} \setminus {\mathbb{V}}^{+})^{\sim }$ is free. It is 
proper because $\mathfrak{G}$ is a discrete subgroup of the Euclidean group 
$\mathrm{E}(2)$ with no accumulation points. \hfill $\square $ \medskip 

\noindent \textbf{Proposition A2.9} $\big( \widehat{\mathfrak{G}} = \langle R \setrule \, R^5 =e \rangle \big) \ltimes 
\mathcal{T}$ is a subgroup of $\mathfrak{G}$, which acts freely, transitively, and properly on 
$(\C \setminus {\mathbb{V}}^{+})^{\sim}$. \medskip 

\noindent \textbf{Proof.} Remove $U$ from the proof of lemmas A2.3 and A2.4. Replace 
$\mathfrak{G}$ by $\widehat{\mathfrak{G}}$ in the proof of lemma A2.7. \hfill $\square $ \medskip 

\noindent \textbf{Theorem A2.10} The $\widehat{\mathfrak{G}}$-orbit space 
$(\mathbb{C} \setminus {\mathbb{V}}^{+})^{\sim }/ \widehat{\mathfrak{G}}$ is holomorphically 
diffeomorphic to the $\widehat{G}$-orbit space $(K^{\ast } \setminus  \mathrm{O} )^{\sim }/\widehat{G} = 
S^{\dagger}$.  \medskip 

\noindent \textbf{Proof} The claim follows because the fundamental 
domain of the $\widehat{\mathfrak{G}}$-action on $ \mathbb{C} \setminus {\mathbb{V}}^{+}$ is 
$K^{\ast } \setminus  \mathrm{O} $. Thus, 
$\Pi ( \mathbb{C} \setminus {\mathbb{V}}^{+})$ 
is a fundamental domain of the $\widehat{\mathfrak{G}}$-action on $(\mathbb{C} \setminus 
{\mathbb{V}}^{+})^{\sim }$, which is equal to $\rho (K^{\ast } \setminus  \mathrm{O})$ 
= $(K^{\ast } \setminus  \mathrm{O} )^{\sim }$. Hence the $\widehat{\mathfrak{G}}$-orbit space $(\mathbb{C} \setminus {\mathbb{V}}^{+})^{\sim }/ \widehat{\mathfrak{G}}$ is equal to the $\widehat{G}$-orbit space $S^{\dagger}$. 
So the identity map from $\Pi (\mathbb{C} \setminus {\mathbb{V}}^{+})$ to 
$(K^{\ast }\setminus  \mathrm{O} )^{\sim }$ induces a holomorphic diffeomorphism 
of orbit spaces. \medskip

\vspace{-.15in}Because the group $\widehat{\mathfrak{G}}$ is a discrete subgroup of the $2$-dimensional 
Euclidean group $\mathrm{E}(2)$, the Riemann surface 
$(\mathbb{C} \setminus {\mathbb{V}}^{+})^{\sim} /\widehat{\mathfrak{G}}$ is an \emph{affine} model of the 
affine Riemann surface $S_{\mathrm{reg}}$.

\end{document}